\documentclass[11pt]{article}
\usepackage[T1]{fontenc}
\usepackage[utf8]{inputenc}
\usepackage{lmodern}
\usepackage[margin=1in]{geometry}
\usepackage{amsmath,amssymb,amsthm,mathrsfs,bbm}
\usepackage{enumitem}
\usepackage{microtype}
\usepackage[colorlinks=true,linkcolor=blue,citecolor=blue,urlcolor=blue]{hyperref}
\usepackage[nameinlink,noabbrev]{cleveref}
\numberwithin{equation}{section}
\newtheorem{theorem}{Theorem}[section]
\newtheorem{proposition}[theorem]{Proposition}
\newtheorem{lemma}[theorem]{Lemma}
\newtheorem{corollary}[theorem]{Corollary}
\theoremstyle{definition}
\newtheorem{definition}[theorem]{Definition}
\newtheorem{remark}[theorem]{Remark}
\newtheorem{conjecture}[theorem]{Conjecture}
\newtheorem{question}[theorem]{Question}
\newcommand{\N}{\mathbb N}
\newcommand{\Z}{\mathbb Z}
\newcommand{\T}{\mathbb T}
\newcommand{\R}{\mathbb R}

\newcommand{\E}{\mathbb E}
\newcommand{\Prob}{\mathbb P}
\newcommand{\ee}{\mathrm e}
\newcommand{\distH}{\operatorname{dist}_{\mathrm H}}
\newcommand{\cB}{\mathscr B}

\title{Sharp $s$-nice recurrence without the van der Corput property}
\author{
Sohail Farhangi\\[0.5em]
Beijing Institute of Mathematical Sciences and Applications (BIMSA)\\
Beijing, China\\
\texttt{sohail.farhangi@gmail.com}
}

\date{}

\begin{document}
\maketitle

\begin{abstract}
For every real $s>2$, we construct a set $R_s\subseteq\N$ which is a set of
$t$-nice measurable recurrence exactly for $t\ge s$ and is not a van der Corput set.
The new ingredient of this construction is the density Frankl--R\"odl theorem.
Furthermore, the probability measure $\sigma_s$ which witnesses failure of the van der Corput property has an atom at zero with mass $(s-2)/(2s-2)$.
A Gaussian tail argument proves both sharpness of the recurrence exponent and the fact that this atom is the largest possible for any $s$-nice recurrent set.
The result sharpens the recurrence separations of Bourgain and Mountakis, answers a question of Mountakis about the maximal mass of an atom at zero, and answers a question of Forrest about separating the various degrees of nice recurrence.

In Appendix A we give a self-contained proof of the density Frankl--R\"odl theorem.
In Appendix B we give a simplified construction of the analogue of $R_s$ in $\mathbb{F}_2^\omega$ instead of $\mathbb{N}$.
\end{abstract}

\section{Introduction}\label{sec:introduction}

We use $(X,\cB,\mu,T)$ to denote an invertible probability-preserving system (p.p.s.), which consists of a standard probability space $(X,\mathscr{B},\mu)$ and a bimeasurable invertible map $T:X\rightarrow X$ satisfying $\mu(TA) = \mu(A)$ for all $A \in \mathscr{B}$.
We let $\mathscr{B}^+ := \{A \in \mathscr{B}\ |\ \mu(A) > 0\}$.
We write $\T=\R/\Z$, let $m_\T$ denote Haar probability measure, and use
$\widehat\sigma(n)=\int_\T e^{-2\pi int}\,d\sigma(t)$ for the Fourier
coefficients of a measure $\sigma$ on $\T$.
We will be studying the relationships between the following notions.

\begin{definition}\label{def:MeasurableRecurrence}
Let $R \subseteq \mathbb{N}$.
\begin{enumerate}
    \item $R$ is a \textbf{set of (measurable) recurrence} if for every p.p.s. $(X,\mathscr{B},\mu,T)$ and every $A \in \mathscr{B}^+$ there exists an $r \in R$ for which $\mu(A\cap T^rA) > 0$.
    
    \item $R$ is a \textbf{set of strong recurrence} if for every p.p.s. $(X,\mathscr{B},\mu,T)$ and every $A \in \mathscr{B}^+$ we have $\displaystyle\limsup_{\underset{r \in R}{r\rightarrow\infty}}\mu(A\cap T^rA) > 0$.
    
    \item Given $s \ge 2$, $R$ is a \textbf{set of $s$-nice recurrence}\footnote{Forrest's original definition of $s$-nice recurrence is different from the one used here, but he proves that it is equivalent to our definition in \cite[Proposition 2.7]{Forrest}.} if for every p.p.s. $(X,\mathscr{B},\mu,T)$ and every $A \in \mathscr{B}^+$ we have $\displaystyle\limsup_{\underset{r \in R}{r\rightarrow\infty}}\mu(A\cap T^rA) \ge \mu(A)^s$.

    \item $R$ is a \textbf{van der Corput (vdC) set}\footnote{While we are not using the original definition of vdC set introduced by Kamae and Mend\`es France \cite{KamaeMendesFrance1978}, our definition is equivalent.
    For a long list of equivalent formulations of vdC sets, the reader is referred to \cite{FarhangiTuckerDrob2026}.} if for every probability measure $\sigma$ on $\T$ satisfying $\widehat{\sigma}(r) = 0$ for every $r \in R$, we have $\sigma(\{0\}) = 0$.
\end{enumerate}
\end{definition}

Khintchine's recurrence theorem \cite{Khintchine1935} tells us that $\mathbb{N}$ is a set of \textbf{nice recurrence}, i.e., a set of $2$-nice recurrence.
It is clear that every set of strong recurrence is a set of recurrence, and that a set of $s$-nice recurrence for some $s \ge 2$ is a set of strong recurrence as well as a set of $t$-nice recurrence for all $t \ge s$.
Bergelson \cite{Bergelson1985} asked if every set of recurrence is a set of strong recurrence, and Forrest \cite{Forrest,Forrest1991} showed that this is not the case (see also \cite{McCutcheon1995}).
Forrest also showed that strong recurrence need not imply nice recurrence by constructing a set of $6.3$-nice recurrence which is not nicely recurrent \cite[Theorem~4.16]{Forrest}.
More generally, for each $t\ge2$, Forrest constructed an $s > t$ and a set of $s$-nice recurrence which is not $t$-nicely recurrent \cite[Corollary~4.17]{Forrest}.
Forrest asked whether, for every $\eta>0$, there exists a $(2+\eta)$-nicely recurrent set that is not $2$-nicely recurrent \cite[p.~156, Question~1]{Forrest}.
Our main result is Theorem \ref{thm:main}, in which we construct for every $s > 2$ a set $R_s$ of $s$-nice recurrence that fails to be $t$-nicely recurrent for every $t < s$, thereby giving a positive answer to the question of Forrest.

Our construction shares some features with Forrest's approach, but the underlying quantitative estimates are different.
Forrest works with finite Boolean groups and constructs sets of allowed differences by imposing a lower bound on Hamming weight.
He then obtains recurrence estimates using combinatorial counting arguments involving binomial coefficients and the binary entropy function.
In contrast, our construction uses differences corresponding to a single prescribed Hamming distance, and we apply the density Frankl--R\"odl theorem to obtain quantitative lower bounds for the proportion of pairs at that distance within an arbitrary subset of the Boolean cube.
The resulting estimates allow us to prescribe the exact exponent $s>2$, rather than merely obtaining a separation between two exponents.
We nevertheless make use of Forrest's Cartesian-power argument \cite[Lemma~2.4 and Proposition~2.7]{Forrest} to eliminate the multiplicative constant appearing in our recurrence estimates.

It is a consequence of the spectral theorem that every vdC set is a set of measurable recurrence, and Ruzsa \cite{Ruzsa1982,Ruzsa1984} asked whether or not the converse holds.
Bourgain \cite{Bourgain} was the first to answer Ruzsa's question by constructing a set of measurable recurrence which is not a vdC set.
Mountakis \cite[Theorem~1.6]{Mountakis} refined Bourgain's construction to produce a set of strong measurable recurrence that is not a vdC set, thereby answering a question of Bergelson and Lesigne \cite{BL}.
We strengthen both of these results by showing that for every $s > 2$, the set $R_s$ of $s$-nice recurrence mentioned above is not a vdC set.
To see that our result is sharp, we recall that recently the author and Rodr\'iguez Mart\'in \cite[Theorem 2.2]{FRM2026} showed that every set of $2$-nice recurrence is a nice vdC set.
In other words, if $R$ is a set of $2$-nice recurrence, then for every probability measure $\sigma$ on $\T$ we have $\displaystyle\limsup_{\underset{r \in R}{r\rightarrow\infty}}|\widehat{\sigma}(r)| \ge \sigma(\{0\})$.

The spirit of our construction is similar to that of Bourgain and Mountakis \cite{Bourgain,Mountakis}.
In all three constructions, separated Riesz products are used to reduce to a finite combinatorial problem involving coordinatewise differences: Bourgain's construction uses the Chinese remainder theorem (CRT), while Mountakis works with base-$Q$ expansions and the finite product group $\mathbb Z_Q^P=(\mathbb Z/Q\mathbb Z)^P$.
In our construction, this reduction takes the form of pairs of vertices in a Boolean cube at a prescribed Hamming distance, with their differences encoded as integers using the Chinese remainder theorem (see also Remark~\ref{rem:TransferFromTorsion}).
Rather than relying on problem-specific combinatorial arguments, we invoke the density Frankl--R\"odl theorem, which provides a quantitative lower bound for the proportion of pairs in a dense subset of the cube that lie at a prescribed Hamming distance.
It is interesting to note that our use of the density Frankl--R\"odl theorem to construct a set of measurable recurrence that is not a vdC set is analogous to the use of Kneser's conjecture \cite{Lovasz1978} to construct a set of topological recurrence that is not a set of measurable recurrence \cite{Griesmer2025Topological,Kriz1987,McCutcheon1995}.

\begin{theorem}\label{thm:main}
For every real $s>2$, there exists a set $R_s$ that is a set of $s$-nice recurrence, but fails to be a set of $t$-nice recurrence for every $t < s$.
Furthermore, there is a symmetric probability measure $\sigma_s$ on $\T$ such that $\widehat\sigma_s(r)=0$ for every $r \in R_s$, while
\begin{equation}\label{eq:main-spectral}
 \sigma_s(\{0\})=\beta_s:= \frac{s-2}{2(s-1)} = \frac{1}{2}-\frac{1}{2(s-1)}.
\end{equation}
Consequently, $R_s$ is not a vdC set.
\end{theorem}

To show that the value of $\sigma_s(\{0\})$ obtained in Equation \eqref{eq:main-spectral} is optimal, we prove the following result.

\begin{theorem}\label{thm:OptimalAtom}
    Let $s > 2$ and let $R$ be a set of $s$-nice recurrence.
    If $\sigma$ is a probability measure on $\T$ for which $\widehat{\sigma}(r) = 0$ for all $r \in R$, then $\sigma(\{0\}) \le \beta_s$.
\end{theorem}

Mountakis' construction shows that for each $\beta \in (0,\frac{1}{3})$ there exists a set of strong recurrence $R$ and a probability measure $\sigma$ on $\T$ for which $\widehat{\sigma}(r) = 0$ for all $r \in R$ and $\sigma(\{0\}) > \beta$.
He asked if the analogous result holds for all $\beta \in (0,\frac{1}{2})$ (cf. \cite[Question~3.5]{Mountakis}). 
We see that Theorem \ref{thm:main} gives a positive answer to this question of Mountakis.

\medskip
\noindent\textbf{Organization of the paper.}
In Section~\ref{sec:recurrence}, we construct the sets $R_s$ using the Chinese remainder theorem and the density Frankl--R\"odl theorem, and prove that they are sets of $s$-nice recurrence.
In Section~\ref{sec:spectral}, we construct the probability measures $\sigma_s$ witnessing the failure of the van der Corput property by means of positive trigonometric polynomials and a separated Riesz product.
In Section~\ref{sec:sharpness}, we use Gaussian processes and tail estimates to prove that the recurrence exponents of the sets $R_s$ and the masses of the atoms $\sigma_s(\{0\})$ are optimal.
In Section~\ref{sec:Questions}, we briefly examine a modification of our construction that allows the sets $R_s$ appearing in Theorem~\ref{thm:main} to form a nested family, and we pose a question and a conjecture for future work.

In Appendix~\ref{app:density-fr}, we give a self-contained proof of the density Frankl--R\"odl theorem, which is the main combinatorial input to our construction.
This theorem is due to Benabbas, Hatami, and Magen; it is stated in \cite[Theorem~5.2]{KOTZ}, and a proof appears in Benabbas's Ph.D.\ thesis \cite[Theorem~7.1 and Section~7.2]{BenabbasThesis}.
We include our proof for completeness and to make the precise formulation needed for our applications readily accessible.

In Appendix~\ref{app:boolean-recurrence} we prove an analogue of Theorem~\ref{thm:main} in $\mathbb{F}_2^\omega$ instead of $\mathbb{Z}$. 
The use of the density Frankl--R\"odl theorem to show $s$-nice recurrence, and the obstruction to the vdC property are both easier to see in this setting.
In Remark~\ref{rem:TransferFromTorsion} we discuss how the CRT blocks construction can be viewed as a tool that transfers the natural construction in $\mathbb{F}_2^\omega$ to $\mathbb{Z}$.

\medskip
\noindent\textbf{AI statement.}
All proofs presented in this paper were originally obtained by GPT-6 Pro, a large language model developed by OpenAI, which also produced the initial draft of the manuscript.
The author subsequently verified the proofs for mathematical correctness and substantially rewrote the manuscript to improve its exposition, clarity, and organization.
The author takes full responsibility for the mathematical content and correctness of the final manuscript.

\medskip
\noindent\textbf{Acknowledgements.} The author would like to thank Yuval Peres for helpful comments on an earlier version of this work.

\section[The construction of \texorpdfstring{$R_s$}{R sub s} and s-nice recurrence]{The construction of \texorpdfstring{$R_s$}{R sub s} and $s$-nice recurrence}
\label{sec:recurrence}

Fix a real number $s>2$ for the rest of the paper.
We construct finite sets $D_j$ and choose integer scales $H_j$ so that $R_s=\bigcup_j H_jD_j$.
The proof that $R_s$ is $s$-nicely recurrent permits arbitrarily large gaps between successive scales $H_j$ and $H_{j+1}$; we use this freedom to separate the Fourier supports in Section~\ref{sec:spectral} to ensure that $R_s$ is not a vdC set.

\subsection{The density Frankl--R\"odl theorem}

For $n\ge1$, give $\{0,1\}^n$ normalized counting measure and write
$\distH(u,v)$ for Hamming distance.
More concretely, we have $\distH(u,v) = \sum_{i = 1}^n|u_i-v_i| = |\{1 \le i \le n\ :\ u_i\neq v_i\}|$.

\begin{theorem}[density Frankl--R\"odl theorem]\label{thm:BHM}
Fix $0<\gamma<1/2$.
There are constants $C_\gamma>0$ and $N_\gamma$ such that the following holds.
Suppose that $n\ge N_\gamma$ and $d=(1-\gamma)n$ is an even integer.
For every $B\subseteq\{0,1\}^n$, writing $a=|B|/2^n$, we have
\begin{equation}\label{eq:BHM}
 \frac{\#\{(x,y)\in B\times B:\operatorname{dist}_{\mathrm H}(x,y)=d\}}
      {2^n\binom nd}
 \ge 2(a/2)^{1/\gamma}
      -C_\gamma\frac{\log^2(n+1)}{n}.
\end{equation}
In particular, the error tends to zero uniformly over all $B$.
The fraction on the left is the probability that both endpoints belong to $B$ when an ordered pair at distance $d$ is chosen uniformly.
\end{theorem}
If $0<\gamma<1/2$ and $(1-\gamma)n$ is an even integer, let $\mathcal E_{n,\gamma}$ be the set of ordered pairs $(u,v)$ at distance $(1-\gamma)n$, with uniform probability measure.
Choose a strictly increasing sequence of rational numbers
\begin{equation}\label{eq:gamma-choice}
 0<\gamma_j<\frac1s,
 \qquad \gamma_j\longrightarrow\frac1s,
 \qquad p_j=\frac1{\gamma_j}>s,
 \qquad c_j=2^{1-p_j}.
\end{equation}
By Theorem~\ref{thm:BHM}, we may choose strictly increasing dimensions $n_j$ such
that
\begin{equation}\label{eq:dimension-choice}
 d_j=(1-\gamma_j)n_j\text{ is an even integer},
 \qquad \varepsilon_j=2^{-j},
\end{equation}
and every $B\subseteq\{0,1\}^{n_j}$ of density $a$ satisfies
\begin{equation}\label{eq:cube-density}
 2^{-n_j}\binom{n_j}{d_j}^{-1}\sum_{(x,y)\in\mathcal E_{n_j,\gamma_j}}\mathbbm{1}_{B\times B}(x,y)
 \ge c_j a^{p_j}-\varepsilon_j.
\end{equation}

\subsection{Constructing \texorpdfstring{$R_s$}{R sub s} with CRT blocks}\label{ssec:ConstructingR_s}

Set
\begin{equation}\label{eq:constants}
 \alpha=1-\frac2s\in(0,1).
\end{equation}
For each $j$, choose $n_j$ distinct odd primes $q_{j,1},\ldots,q_{j,n_j}$
so large that
\begin{equation}\label{eq:prime-condition}
 W_j:=\sum_{i=1}^{n_j}\frac1{1+\cos(\pi/q_{j,i})}
 <\frac{d_j}{1+\alpha}.
\end{equation}
This is possible because $W_j\to n_j/2$ as the primes tend to infinity,
whereas
\begin{equation*}
 \frac{d_j}{1+\alpha}
 =\frac{(1-\gamma_j)n_j}{2(1-1/s)}>\frac{n_j}{2}.
\end{equation*}

Set $Q_j=\prod_{i=1}^{n_j}q_{j,i}$ and let
$\Phi_j:\{0,1\}^{n_j}\to\{0,\ldots,Q_j-1\}$ be the injective map
specified by
\begin{equation*}
 \Phi_j(v)\equiv v_i\pmod{q_{j,i}}
 \qquad(1\le i\le n_j).
\end{equation*}
Define
\begin{equation*}
\begin{split}
 D_j=\{r\in\{1,\ldots,Q_j-1\}:\;&r \equiv 0,1,\text{ or }-1\bmod q_{j,i}
                         \text{ for every }i,\#\{i:r\not\equiv0\bmod{q_{j,i}}\}=d_j\}.
\end{split}
\end{equation*}
These sets are nonempty and $Q_j-D_j=D_j$.  If
$\distH(u,v)=d_j$, then
\begin{equation}\label{eq:edge-difference}
 |\Phi_j(u)-\Phi_j(v)|\in D_j.
\end{equation}
In particular, $D_j$ contains the actual positive integer difference,
including the representative that can arise from wraparound modulo
$Q_j$.

Lemma~\ref{lem:local-polynomial} constructs, from these local data alone,
a positive trigonometric polynomial $P_j$ of degree at most $L_j$ with
$\widehat P_j(0)=1$ and $\widehat P_j(r)=-\alpha$ on $D_j$.
Its construction is independent of the scales $H_j$.
Consequently, we may choose the scales $H_j$ recursively so that
\begin{equation}\label{eq:scale-separation}
 H_1=1,\qquad
 H_{j+1}>3\sum_{i=1}^j L_iH_i.
\end{equation}
We mention once again that the proof that $R_s = \bigcup_{j = 1}^\infty H_jD_j$ is $s$-nicely recurrent works for every choice of $H_j$.
We only require a growth condition on the $H_j$ to show that $R_s$ is not a vdC set, and that $R_s$ fails to be $t$-nicely recurrent for every $t < s$.

\subsection{s-nice recurrence from orbit averaging and Cartesian powers}

The fact that $R_s$ is a set of $s$-nice recurrence follows immediately from the following result.

\begin{proposition}\label{prop:block-recurrence}
For a p.p.s. $(X,\mathscr{B},\mu,T)$ and a measurable set $A \in \mathscr{B}$, write
\begin{equation*}
 M_j(A)=\max_{r\in H_jD_j}\mu(A\cap T^{-r}A).
\end{equation*}
For every positive integer $k$, we have
\begin{equation}\label{eq:tensor-bound}
 M_j(A)^k\ge c_j\mu(A)^{kp_j}-\varepsilon_j.
\end{equation}
Consequently,
\begin{equation}\label{eq:block-recurrence}
 \liminf_{j\to\infty}M_j(A)\ge\mu(A)^s.
\end{equation}
\end{proposition}

\begin{proof}
For $x\in X$, let
\begin{equation*}
 B_j(x)=\left\{v\in\{0,1\}^{n_j}:T^{-H_j\Phi_j(v)}x\in A\right\},
 \qquad a_j(x)=\frac{|B_j(x)|}{2^{n_j}}.
\end{equation*}
We see that

\begin{equation*}
    \int a_jd\mu = 2^{-n_j}\sum_{v \in \{0,1\}^{n_j}}\int_X\mathbbm{1}_A\left(T^{-H_j\Phi_j(v)}x\right)d\mu(x) = \mu(A).
\end{equation*}
Recalling that for every $(u,v) \in \mathcal{E}_{n_j,\gamma_j}$ we have $|\Phi_j(u)-\Phi_j(v)| \in D_j$, we use Theorem~\ref{thm:BHM} and Jensen's inequality to see that
\begin{alignat*}{2}
 M_j(A)
 &\ge 2^{-n_j}\binom{n_j}{d_j}^{-1}\sum_{(u,v) \in \mathcal{E}_{n_j,\gamma_j}}\mu\left(T^{H_j\Phi_j(u)}A\cap T^{H_j\Phi_j(v)}A\right)\\
 &= \int_X2^{-n_j}\binom{n_j}{d_j}^{-1}\sum_{(u,v) \in \mathcal{E}_{n_j,\gamma_j}}\mathbbm{1}_A\left(T^{-H_j\Phi_j(u)}x\right)\mathbbm{1}_A\left(T^{-H_j\Phi_j(v)}x\right)d\mu(x)\\
 &= 2^{-n_j}\binom{n_j}{d_j}^{-1}\int_X\sum_{(u,v)\in\mathcal E_{n_j,\gamma_j}}\mathbbm{1}_{B_j(x)\times B_j(x)}(u,v)\,d\mu(x) \ge c_j\int_X a_j(x)^{p_j}\,d\mu(x)-\varepsilon_j\\
 &\ge c_j\mu(A)^{p_j}-\varepsilon_j.
\end{alignat*}
Apply this inequality to the Cartesian-power system $(X^k,\mu^{\otimes k},T^{\times k})$ and the set $A^k$.
At every fixed time, its return correlation is the $k$th power of the original one.
Thus $M_j(A^k)=M_j(A)^k$, proving \eqref{eq:tensor-bound}.

For fixed $k$, let $j\to\infty$.
Since $p_j\to s$, $c_j\to2^{1-s}$, and $\varepsilon_j\to0$, we obtain
\begin{equation*}
 \liminf_j M_j(A)\ge 2^{(1-s)/k}\mu(A)^s.
\end{equation*}
Now we let $k\to\infty$ to prove \eqref{eq:block-recurrence}.
\end{proof}

We mention that the Cartesian-power removal of a fixed positive multiplicative constant is the argument of Forrest \cite[Lemma~2.4 and Proposition~2.7]{Forrest}.

\section{The spectral construction and failure of the vdC property}
\label{sec:spectral}

Put $\ee(t)=e^{2\pi it}$.
We construct a symmetric probability measure $\nu$ with $\widehat\nu(r)=-\alpha$ for every $r\in R_s$, where $\alpha=1-2/s$.
Taking a convex combination with an atom at zero then gives the spectral witness.

\subsection{A finite Fourier identity}

\begin{lemma}\label{lem:local-polynomial}
For each block $D_j$ there is a real, even trigonometric polynomial $P_j$
of degree at most $L_j$, where $L_j\ge Q_j$, such that
\begin{equation}\label{eq:local-polynomial}
 P_j>0,\qquad \widehat P_j(0)=1,\qquad
 \widehat P_j(r)=-\alpha\quad(r\in D_j).
\end{equation}
\end{lemma}

\begin{proof}
We suppress the block index by writing $n=n_j$, $d=d_j$, $q_i=q_{j,i}$, $D = D_j$, and $Q = Q_j$, then put
\begin{equation*}
 t_i=\frac{q_i-1}{2q_i},\qquad
 w_i=\frac1{1+\cos(\pi/q_i)},\qquad W=\sum_{i=1}^n w_i,
 \qquad
 \eta_0=\frac12\sum_{i=1}^n w_i(\delta_{t_i}+\delta_{-t_i}).
\end{equation*}
If $r\equiv0\pmod{q_i}$, the $i$th summand contributes $w_i$ to $\widehat\eta_0(r)$.
If $r\equiv\pm1\pmod{q_i}$, it contributes $-w_i\cos(\pi/q_i)=w_i-1$.
Every $r\in D$ has exactly $d$ nonzero CRT coordinates, so
\begin{equation}\label{eq:finite-Fourier-identity}
 \widehat\eta_0(r)=W-d\qquad(r\in D).
\end{equation}
By \eqref{eq:prime-condition}, $W<d/(1+\alpha)$.  Hence
$b:=\alpha W/(d-W)<1$, and
\begin{equation}\label{eq:local-positive-measure}
 \lambda=\frac{\alpha}{d-W}\eta_0+(1-b)m_\T
\end{equation}
is a symmetric probability measure with
$\widehat\lambda(r)=-\alpha$ on $D$ and a positive Haar component.

Let $F_L$ be the Fej\'er kernel, so that
\begin{equation*}
 F_L(t)=\frac1{L+1}\left|\sum_{u=0}^L\ee(ut)\right|^2,
 \qquad
 \widehat F_L(a)=\left(1-\frac{|a|}{L+1}\right)_+.
\end{equation*}
For $L\ge Q$, correct the finitely many prescribed coefficients by setting
\begin{equation}\label{eq:polynomial-correction}
 P=\lambda*F_L+
 \sum_{r\in D\cup(-D)}
   \frac{|r|}{L+1}\widehat\lambda(r)\ee(rt).
\end{equation}
The correction is real and even, has mean zero, and its uniform norm is at most $2\alpha\sum_{r\in D}r/(L+1)$, which tends to zero.
Since $\lambda*F_L\ge1-b>0$, a sufficiently large $L$ makes $P$ strictly positive.
Its mean is one and its coefficients on $D$ are exactly $-\alpha$.
Take $L_j$ to be such a degree bound, enlarged if necessary.
\end{proof}

\subsection{A separated Riesz product}

For $m\ge1$, set
\begin{equation*}
 \Pi_m(t)=\prod_{j=1}^mP_j(H_jt),\qquad
 S_m=\sum_{j=1}^m L_jH_j,\qquad S_0=0.
\end{equation*}
Here a Riesz product means the weak limit of these positive polynomial products.
The separation \eqref{eq:scale-separation} makes their Fourier coefficients explicit.

\begin{lemma}\label{lem:dissociation}
The probability measures $\Pi_m\,dm_\T$ converge weakly to a symmetric
probability measure $\nu$.  Every integer has at most one finite
representation
\begin{equation*}
 a=\sum_j e_jH_j,\qquad e_j\in\Z,\qquad |e_j|\le L_j.
\end{equation*}
For this representation,
\begin{equation}\label{eq:nu-coefficients}
 \widehat\nu(a)=\prod_j\widehat P_j(e_j),
\end{equation}
and $\widehat\nu(a)=0$ if no such representation exists.  In particular,
\begin{equation}\label{eq:negative-correlations}
 \widehat\nu(H_jr)=-\alpha\quad(r\in D_j),\qquad
 \widehat\nu(r)=-\alpha\quad(r\in R_s).
\end{equation}
\end{lemma}

\begin{proof}
For every $k\ge1$, \eqref{eq:scale-separation} gives
$S_{k-1}<H_k/3$.  If two representations differ, let $k$ be the largest
index of a nonzero digit difference.  That digit contributes at least
$H_k$ in absolute value, while all smaller digits contribute at most
$2S_{k-1}<2H_k/3$, a contradiction.

Expanding $\Pi_m$ therefore gives the asserted product formula for its
coefficients.  In particular, its constant coefficient is one, so
$\Pi_m\ge0$ is a probability density.  A nonzero frequency whose largest
digit has index $k$ has absolute value at least
$H_k-S_{k-1}>2H_k/3$.  Thus every fixed coefficient stabilizes as
$m\to\infty$.  Weak compactness of probability measures on $\T$ and
uniqueness from Fourier coefficients give the asserted limit and
\eqref{eq:nu-coefficients}.  Evenness of the $P_j$ gives symmetry of
$\nu$.  The one-digit representation of $H_jr$ proves
\eqref{eq:negative-correlations}.
\end{proof}

\subsection{The measures \texorpdfstring{$\sigma_s$}{sigma sub s}}\label{ssec:TheMeasures}

Define
\begin{equation}\label{eq:spectral-witness}
 \sigma_s=\frac{\alpha\delta_0+\nu}{1+\alpha}.
\end{equation}
This is a symmetric probability measure and, for every $r\in R_s$,
\begin{equation}\label{eq:witness-atom-lower}
 \widehat\sigma_s(r)=\frac{\alpha+\widehat\nu(r)}{1+\alpha}=0,
 \qquad
 \sigma_s(\{0\})\ge\frac{\alpha}{1+\alpha}
 =\frac{s-2}{2(s-1)}>0.
\end{equation}
Theorem~\ref{thm:OptimalAtom}, which we prove below, together with the $s$-nice recurrence already proved, forces equality in \eqref{eq:witness-atom-lower}.

\section{Sharpness of the recurrence exponent and the spectral atom}
\label{sec:sharpness}

A single Gaussian estimate proves both optimality statements.
We first state it for an arbitrary real covariance sequence, and then apply it to the measure constructed above and to an arbitrary spectral obstruction.

\subsection{A Gaussian tail estimate}

\begin{proposition}[Gaussian tail obstruction]
\label{prop:gaussian-obstruction}
Let $E\subseteq\N$ be nonempty, and let $c\colon\Z\to\R$ be a real, even, positive-definite sequence satisfying
\begin{equation}
 c(0)=1,\qquad c(r)=-a\quad(r\in E),\qquad 0\le a<1.
 \label{eq:gaussian-obstruction-covariance}
\end{equation}
Put $q=2/(1-a)$.
There is a probability-preserving system $(X,\mathscr{B},\mu,T)$ with measurable sets $A_u$, $u\ge1$, such that
\begin{equation}
 \mu(A_u) \ge \frac{e^{-3/2}}{\sqrt{2\pi}}u^{-1}\exp\left(-\frac{u^2}{2}\right),
 \qquad \sup_{r\in E}\mu(A_u\cap T^{-r}A_u)
 \le \exp\left(-\frac{q u^2}{2}\right).
 \label{eq:gaussian-obstruction-joint-tail}
\end{equation}
Consequently, for every real  $2 \le t < q$,
\begin{equation}
 \lim_{u\to\infty}\displaystyle\sup_{r\in E}
 \frac{\mu(A_u\cap T^{-r}A_u)}
      {\mu(A_u)^t}=0,
 \label{eq:gaussian-obstruction-ratio}
\end{equation}
hence $E$ is not a set of $t$-nice recurrence for any $2 \le t < q$.
\end{proposition}

\begin{proof}
We briefly recall the Gaussian measure space construction here, and we refer the reader to \cite[Chapter 8.2]{CornfeldFominSinai1982} and \cite[Chapter 3.11]{Glasner2003} for more detailed treatments.
 For every finite set of integer indices $n_1,\ldots,n_k$, positive definiteness makes $(c(n_i-n_j))_{i,j=1}^k$ a positive semidefinite real matrix.
 The centered Gaussian distributions with these covariance matrices form a consistent family.
The Kolmogorov extension theorem therefore gives a probability measure $\mu$ on $X=\R^\Z$ under which the coordinate functions $(Y_n)_{n\in\Z}$ are jointly centered Gaussian and satisfy
\begin{equation}
 \E(Y_nY_m)=c(n-m),\qquad \E(Y_n^2)=1.
 \label{eq:gaussian-obstruction-realization}
\end{equation}
Simultaneously translating the indices leaves the covariance matrices unchanged, so the law is invariant under the shift $(Tx)_n=x_{n+1}$.
In particular, $Y_n=Y_0\circ T^n$.

For $u \ge 1$ let $A_u = \{Y_0 > u\}$.
For each $r\in E$, the centered Gaussian variable $Y_0+Y_r$ has variance $2(1-a)=4/q$.
For any real random variable $Z$ and any $\lambda > 0$, Markov's inequality applied to $e^{\lambda Z}$ gives $\Prob(Z>x)\le e^{-\lambda x}\E e^{\lambda Z}$.
If $Z$ is a centered Gaussian variable with variance $v$, then $\E e^{\lambda Z}=e^{v\lambda^2/2}$ (cf. \cite[p.~22]{BLM}).
Taking $Z=Y_0+Y_r$ gives
\begin{alignat*}{2}
 \mu(A_u\cap T^{-r}A_u)
 &\le \mu(\{Y_0+Y_r >2u\}) \le \exp(-2\lambda u)\E\exp\bigl(\lambda(Y_0+Y_r)\bigr) =\exp\left(-2\lambda u+\frac{2\lambda^2}{q}\right).
\end{alignat*}
Taking $\lambda=qu/2$ proves
\eqref{eq:gaussian-obstruction-joint-tail}, uniformly in $r\in E$.

For $u\ge1$, integration of the standard normal density over $[u,u+1/u]$ yields
\begin{alignat*}{2}
 \mu(A_u)
 &\ge\frac1{\sqrt{2\pi}\,u}
       \exp\left(-\frac{(u+1/u)^2}{2}\right) \ge \frac{e^{-3/2}}{\sqrt{2\pi}}u^{-1}\exp\left(-\frac{u^2}{2}\right).
\end{alignat*}
\end{proof}

\subsection{The optimal mass at zero}

\begin{proof}[Proof of Theorem \ref{thm:OptimalAtom}]
Write $\beta=\sigma(\{0\})$.
For any $r\in R$, the equation $\widehat\sigma(r)=0$ implies
\begin{equation*}
 \beta
 =\left|\int_{\T\setminus\{0\}}e^{-2\pi i rt}\,d\sigma(t)\right|
 \le1-\beta,
\end{equation*}
hence $\beta\le1/2$.
In fact, we claim that $\beta<1/2$.
To see this, let $r_1 \in R$ be arbitrary and note that we can only have $\beta = 1/2$ and $\widehat{\sigma}(r_1) = 0$ if $\sigma-\frac12\delta_0$ is supported on the $2r_1^{\text{th}}$ roots of unity.
Since $R$ is a set of recurrence, it contains a multiple of every integer, so let $r_2 \in R$ be a multiple of $2r_1$.
We see that $\widehat{\sigma}(r_2) = 1$, which contradicts our assumption that $\widehat{\sigma}(r_2) = 0$ and proves the claim.

Now replace $\sigma$ by its symmetrization $(\sigma+\check\sigma)/2$, where $\check\sigma(B)=\sigma(-B)$.
This preserves its mass at zero and its vanishing Fourier coefficients on $R$.
The measure
\begin{equation*}
 \tau=\frac{\sigma-\beta\delta_0}{1-\beta}
\end{equation*}
is now a symmetric probability measure, and its Fourier coefficients form a real, even, positive-definite sequence with
\begin{equation}
 \widehat\tau(0)=1,\qquad
 \widehat\tau(r)=-\frac{\beta}{1-\beta}\quad(r\in R).
 \label{eq:optimal-atom-remainder}
\end{equation}

Set
\begin{equation*}
 a=\frac{\beta}{1-\beta},\qquad
 q=\frac2{1-a}=\frac{2(1-\beta)}{1-2\beta}.
\end{equation*}
If $\beta>(s-2)/[2(s-1)]$, then $q>s$, so we can apply Proposition~\ref{prop:gaussian-obstruction} to $c(n)=\widehat\tau(n)$ with $t=s$ to see that $R$ is not a set of $s$-nice recurrence, a contradiction.
\end{proof}

\subsection{The exact recurrence exponent}

\begin{corollary}\label{cor:sharpness}
The set $R_s$ constructed in Section \ref{ssec:ConstructingR_s} fails to be a set of $t$-nice recurrence for every $t < s$.
\end{corollary}

\begin{proof}
    Let $2 < t < s$ be arbitrary and let $\sigma_s$ be as in Section \ref{ssec:TheMeasures}.
    We showed in Section \ref{ssec:TheMeasures} that

    \begin{equation}
        \sigma_s(\{0\}) \ge \frac{1}{2}-\frac{1}{2(s-1)} > \frac{1}{2}-\frac{1}{2(t-1)},
    \end{equation}
    so Theorem \ref{thm:OptimalAtom} tells us that $R_s$ is not a set of $t$-nice recurrence.
\end{proof}

\section{Questions and future work}\label{sec:Questions}

We observe that it is possible to strengthen Theorem~\ref{thm:main} by requiring the sets $R_s$ to form a nested family, in the sense that $R_s\subsetneq R_t$ whenever $2 < t < s$.
To obtain this strengthening, one can replace the fixed parameter $s$ in Section~\ref{sec:recurrence} with a sequence $(s_j)_{j\ge1}$ enumerating all rational numbers greater than $2$.
For each $j$, put $p_j=s_j+1/j$, $\gamma_j=1/p_j$, and $\alpha_j'=1-2/s_j$, and construct the corresponding CRT block $D_j$ and positive trigonometric polynomial $P_j$ using these parameters.
At each stage, choose the dimension sufficiently large that the cube-density estimate holds with error $\varepsilon_j=2^{-j}$.
The strict inequality $p_j>s_j$ ensures that the prime-selection condition \eqref{eq:prime-condition} can still be satisfied.
Choose the scales $H_j$ to satisfy \eqref{eq:scale-separation}, independently of $s$, and define
\[
R_s=\bigcup_{\substack{s_j\ge s}}H_jD_j
\qquad(s>2).
\]
The desired nesting follows immediately from this definition, and the inclusions are strict because the blocks are disjoint and the rational numbers are dense.
For each fixed $s>2$, the proof of $s$-nice recurrence can be recovered by applying the orbit-averaging and Cartesian-power arguments along a subsequence of indices for which $s_j>s$ and $s_j\to s$.

The spectral construction requires one additional modification.
For each $s>2$, let $\alpha_s=1-2/s$ and replace $P_j$ with
\[
P_{j,s}=
\begin{cases}
\displaystyle
1+\frac{\alpha_s}{\alpha_j'}(P_j-1),
& s_j\ge s,\\[6pt]
1,&s_j<s.
\end{cases}
\]
Since $\alpha_s\le\alpha_j'$ whenever $s_j\ge s$, each $P_{j,s}$ is a convex combination of $P_j$ and the constant polynomial $1$.
Consequently, these polynomials remain strictly positive, have constant Fourier coefficient one, and satisfy $\widehat P_{j,s}(r)=-\alpha_s$ for every $r\in D_j$ whenever $s_j\ge s$.
Their degrees do not increase, so the separated Riesz-product argument of Section~\ref{sec:spectral} continues to apply.
This produces a symmetric probability measure $\sigma_s$ whose Fourier coefficients vanish on $R_s$ and whose atom at zero has the optimal mass $\beta_s$.
Theorem~\ref{thm:OptimalAtom} then gives the sharpness of the recurrence exponent exactly as before.

We believe that the sets presented by Theorem \ref{thm:main} are plentiful, which leads us to the following conjecture.

\begin{conjecture}\label{conj:Hereditary}
    Suppose that $R\subseteq\mathbb N$ is a set of $2$-nice recurrence.
    There exists a family $(R_s)_{s>2}$ of subsets of $R$ with the following properties.
    \begin{enumerate}
        \item $R_s$ is a set of $s$-nice recurrence, but fails to be a set of $t$-nice recurrence for every $t < s$.
        \item There exists a probability measure $\sigma_s$ on $\T$ for which $\widehat{\sigma_s}(r) = 0$ for every $r \in R_s$ and $\sigma_s(\{0\}) = \beta_s$.
        \item If $2 < t < s$, then $R_s \subsetneq R_t$.
    \end{enumerate}
\end{conjecture}
Some partial progress towards Conjecture~\ref{conj:Hereditary}
has been obtained using GPT-6.
In particular, the AI-generated manuscript \cite{GPT6RecurrenceExtraction}
contains claimed proofs of the conjecture in several special cases of
interest, including suitable polynomial images of IP-sets,
difference sets, and sets arising from Hartman sequences.\footnote{
A \textbf{difference set} is a set of the form
$\{n_i-n_j\mid i>j\}$ for some increasing sequence
$(n_i)_{i=1}^\infty$ in $\mathbb{N}$.
A sequence $(n_i)_{i=1}^\infty$ is \textbf{Hartman} if for every
$\alpha\in\mathbb{R}$, the sequence
$(n_i\alpha\pmod{1})_{i=1}^\infty$ is uniformly distributed
with respect to Haar measure on the orbit closure of
$\{n\alpha\pmod{1}\mid n\in\mathbb{N}\}$.}
The proofs in \cite{GPT6RecurrenceExtraction} have not been
fully verified by the author, and we mention them here as
possible avenues for future work.

A set $A\subseteq\mathbb{N}$ is an \textbf{IP-set} if there
exists a sequence $(x_n)_{n=1}^\infty$ in $\mathbb{N}$ such that
every nonempty finite sum of distinctly indexed elements of
$(x_n)_{n=1}^\infty$ is contained in $A$.
Similarly, $A$ is an \textbf{$\mathrm{IP}_0$-set} if for every
$N\in\mathbb{N}$ there exists a sequence $(x_{n,N})_{n=1}^N$
in $\mathbb{N}$ such that every nonempty finite sum of distinctly
indexed elements of $(x_{n,N})_{n=1}^N$ is contained in $A$.
A set is an $\mathrm{IP}_r^*$-set if it has nonempty intersection
with every finite-sums set generated by $r$ positive integers,
and is an $\mathrm{IP}_0^*$-set if it has nonempty intersection
with every $\mathrm{IP}_0$-set, or equivalently, if it is an
$\mathrm{IP}_r^*$-set for some $r\in\mathbb{N}$.

A classical source of sets of $2$-nice recurrence is the
polynomial images of IP-sets.
More concretely, Bergelson, Furstenberg, and McCutcheon
\cite[Corollary~2.1]{BergelsonFurstenbergMcCutcheon1996}
showed that if $A$ is an IP-set and
$p\in x\mathbb{Z}[x]\setminus\{0\}$, then $p(A)$ is a set
of $2$-nice recurrence.\footnote{
Throughout this discussion, when considering the recurrence
properties of a set $S\subseteq\mathbb{Z}$, we mean the
corresponding properties of
$\{|n|\mid n\in S\setminus\{0\}\}\subseteq\mathbb{N}$.}

There are also results concerning the recurrence properties of polynomial images of $\mathrm{IP}_0$-sets.
In particular, a special case of the work of Bergelson and Leibman \cite[Theorem~4.9]{BergelsonLeibman2018Cubics} (see also \cite[Theorem~1.16 and Remark~1.17]{BergelsonZelada2026}) shows that for every p.p.s. $(X,\mathscr{B},\mu,T)$, every $B\in\mathscr{B}^+$, every $p\in x\mathbb{Z}[x]\setminus\{0\}$, and every $\varepsilon>0$, the set
\[
R_\varepsilon^p(B) := \{n\in\mathbb{N}\mid \mu(B\cap T^{-p(n)}B)>\mu(B)^2-\varepsilon\}
\]
is an almost $\mathrm{IP}_0^*$-set, meaning that it becomes an $\mathrm{IP}_0^*$-set upon adjoining a suitable set of zero upper Banach density.
In fact, the aforementioned results establish the stronger almost $\mathrm{VIP}_0^*$ property (cf. \cite[Definition~2.2]{BergelsonLeibman2018Cubics}).

Nevertheless, Zelada \cite[Theorem~A]{Zelada2024} showed that for every $p\in x\mathbb{Z}[x]$ with $\deg(p)>1$, there exists an $\mathrm{IP}_0$-set $A_p$ such that $p(A_p)$ is not a set of $2$-nice recurrence.
It remains to be determined whether the exceptional set of zero upper Banach density in the aforementioned results can be dispensed with if one only asks for positive intersections rather than intersections of size approximately $\mu(B)^2$.
This leads us to the following question, which can be shown to be equivalent to Question~1.4 of Zelada \cite{Zelada2024}.

\begin{question}\label{q:BabyPDHJ}
Suppose that $A\subseteq\mathbb{N}$ is an $\mathrm{IP}_0$-set and $p\in x\mathbb{Z}[x]\setminus\{0\}$.
Is $p(A)$ a set of measurable recurrence?
\end{question}

Question~\ref{q:BabyPDHJ} can be viewed as the most basic case of Bergelson's \cite[p.~41]{Bergelson2010} question about whether the $\mathrm{IP}_r^*$-polynomial Szemer\'edi theorem holds.
The connection with the density polynomial Hales--Jewett conjecture (cf. \cite[Conjecture~3]{Gowers2009} and \cite[Conjecture~1.1]{Karam2026}) has recently been made particularly explicit by Bergelson and Zelada \cite[Theorem~A1]{BergelsonZelada2026}.
Their result shows that the density polynomial Hales--Jewett conjecture implies a quantitative form of the $\mathrm{IP}_r^*$-polynomial Szemer\'edi theorem for
commuting transformations.

In connection with the setting of Appendix~\ref{app:boolean-recurrence}, we also mention the recent work of Zelada \cite[Theorem~C1]{Zelada2026Polynomial}.
He establishes an equivalence between a symmetric-difference weakening of the density polynomial Hales--Jewett conjecture and several polynomial recurrence statements for $\mathbb{F}_2^\omega$-actions, including a quantitative finite VIP-recurrence theorem.
These results provide further motivation for studying single polynomial recurrence along $\mathrm{IP}_0$-sets, both in $\mathbb{Z}$ and in more general countable abelian groups.

\appendix
\section{A self-contained proof of the density Frankl--R\"odl theorem}
\label{app:density-fr}

We prove the density Frankl--R\"odl theorem as formulated in Theorem \ref{thm:BHM}.
The underlying density theorem is due to Benabbas, Hatami, and Magen; it is stated in \cite[Theorem~5.2]{KOTZ}, and a proof appears in Benabbas's Ph.D.\ thesis \cite[Theorem~7.1 and Section~7.2]{BenabbasThesis}.
Our appendix gives a self-contained proof with the explicit error bound used in Theorem~\ref{thm:BHM}.
We use the strategy of Benabbas of comparing the average of two consecutive Hamming-distance operators with an independent-noise operator.

In Section A.1 we introduce notation for the basic objects that we will be working with.
In Sections A.2 and A.3 we prove technical inequalities that will be used in Section A.4 to prove Lemma \ref{appfr:comparison}, the key operator comparison that underpins the main proof.
In Section A.5 we put everything together and prove the density Frankl--R\"odl theorem.

We mention here that the original Frankl--R\"odl theorem \cite{FranklRodl1987} starts with a given distance $d$ and bounds from above the size of sets $A \subseteq \{0,1\}^n$ that contain no pair of points with Hamming distance $d$.
The density Frankl--R\"odl theorem starts with a given set $A \subseteq \{0,1\}^n$ and a given distance $d$ and bounds from below the number of pairs in $A$ having Hamming distance $d$.

\subsection{Hamming cubes, Fourier characters, and averaging operators}

Write $\Omega_n=\{0,1\}^n$, with coordinatewise addition modulo two and normalized counting measure.
For real functions on $\Omega_n$, put
\begin{equation*}
 \mathbb E f=2^{-n}\sum_{x\in\Omega_n}f(x),
 \qquad
 \langle f,g\rangle=\mathbb E(fg),
 \qquad
 \|f\|_2=\langle f,f\rangle^{1/2}.
\end{equation*}
For a nonnegative function $f$ and $p>0$, also write
\begin{equation*}
 \|f\|_p=(\mathbb E f^p)^{1/p}.
\end{equation*}
Let $|x|=\sum_{i=1}^n x_i$ and $[n]=\{1,\ldots,n\}$.

For $S\subseteq[n]$, define the Walsh character
\begin{equation*}
 \chi_S(x)=(-1)^{\sum_{i\in S}x_i}.
\end{equation*}
The identity $\chi_S\chi_T=\chi_{S\mathbin\triangle T}$ shows that these characters are orthonormal.
There are $2^n$ of them, so they form a basis for the $2^n$-dimensional space of real functions on $\Omega_n$.
Consequently,
\begin{equation}
\label{appfr:fourier}
 f=\sum_{S\subseteq[n]}\widehat f(S)\chi_S,
 \qquad
 \widehat f(S)=\mathbb E(f\chi_S),
 \qquad
 \|f\|_2^2=\sum_{S\subseteq[n]}\widehat f(S)^2.
\end{equation}

For $0\le d\le n$, define the spherical averaging operator by
\begin{equation*}
 (\mathsf S_df)(x)
   =\binom nd^{-1}\sum_{\substack{z\in\Omega_n\\ |z|=d}}f(x+z).
\end{equation*}
For $-1\le\rho\le1$, define the noise operator by
\begin{equation*}
 (\mathsf T_\rho f)(x)=\mathbb E_Z f(x+Z),
 \qquad
 \mathbb P(Z_i=1)=\frac{1-\rho}{2},
\end{equation*}
where the coordinates of $Z$ are independent.
Thus $\mathsf S_d$ averages over vertices that differ from $x$ in exactly $d$ coordinates, whereas $\mathsf T_\rho$ averages over a random perturbation of the coordinates of $x$.
It can be checked that
\begin{equation}
\label{appfr:eigenvalues}
 \mathsf T_\rho\chi_S=\rho^{|S|}\chi_S,
 \qquad
 \mathsf S_d\chi_S=\lambda_d(|S|)\chi_S,
\end{equation}
where
\begin{equation}
\label{appfr:lambda}
 \lambda_d(k)
 =\binom nd^{-1}
   \sum_{j=0}^k(-1)^j\binom kj\binom{n-k}{d-j}.
\end{equation}
Here and below a binomial coefficient is zero when its lower index is outside its usual range.

In particular, for $f=\mathbf 1_B$,
\begin{equation}
\label{appfr:edge-probability}
 \langle f,\mathsf S_df\rangle
 =\frac{\#\{(x,y)\in B\times B:
                  \operatorname{dist}_{\mathrm H}(x,y)=d\}}
        {2^n\binom nd}.
\end{equation}

Equation \eqref{appfr:edge-probability} shows us that understanding $S_d$ is the key to proving the density Frankl--R\"odl theorem, but Equation \eqref{appfr:eigenvalues} shows us that it is significantly easier to work with $\mathsf T_\rho$ than $\mathsf S_d$.
Most of the work of this proof goes towards proving Lemma \ref{appfr:comparison}, which shows that for $\rho = 1-2d/n$, the operator $\mathsf T_\rho$ is close in operator norm to the average of $\mathsf S_d$ and $\mathsf S_{d+1}$.

\subsection{An elementary reverse-hypercontractive inequality}

We first prove the precise noise inequality that will be needed.

\begin{lemma}
\label{appfr:scalar}
If $0<p\le1$ and $0\le t\le1$, then
\begin{equation}
\label{appfr:scalar-bound}
 \left(\frac{(1+t)^p+(1-t)^p}{2}\right)^{2/p}
 \le 1-(1-p)t^2.
\end{equation}
\end{lemma}

\begin{proof}
The case $p=1$ is immediate.
Suppose $0<p<1$, and put $q=1-p$.
For $0\le t<1$, define
\begin{equation*}
 F(t)=\log\left(\frac{(1+t)^p+(1-t)^p}{2}\right)
       -\frac p2\log(1-qt^2).
\end{equation*}
Differentiating and combining the two terms gives
\begin{equation*}
 F'(t)
 =p\frac{(1+qt)(1+t)^{-q}-(1-qt)(1-t)^{-q}}
 {\bigl((1+t)^p+(1-t)^p\bigr)(1-qt^2)}.
\end{equation*}
The denominator is positive.
To see that the numerator is nonpositive, let
\begin{equation*}
 G(t)=q\log\frac{1+t}{1-t}
      -\log\frac{1+qt}{1-qt}.
\end{equation*}
We have $G(0)=0$ and
\begin{equation*}
 G'(t)
 =2q\left(\frac1{1-t^2}-\frac1{1-q^2t^2}\right)\ge0.
\end{equation*}
Hence
\begin{equation*}
 \left(\frac{1+t}{1-t}\right)^q
 \ge\frac{1+qt}{1-qt},
\end{equation*}
which is exactly the inequality needed to make the numerator of $F'(t)$ nonpositive.
Thus $F(t)\le F(0)=0$.
Exponentiating and raising to the power $2/p$ proves \eqref{appfr:scalar-bound} for $t<1$.
Continuity gives the endpoint $t=1$.
\end{proof}

\begin{lemma}[Reverse hypercontractivity]
\label{appfr:reverse}
Let $-1<\rho<1$, and put $p=1-|\rho|$.
For all nonnegative $f,g:\Omega_n\to\mathbb R$,
\begin{equation}
\label{appfr:reverse-bound}
 \langle f,\mathsf T_\rho g\rangle\ge\|f\|_p\|g\|_p.
\end{equation}
In particular, if $A\subseteq\Omega_n$ has density $b$, then
\begin{equation}
\label{appfr:noise-set-bound}
 \langle\mathbf 1_A,\mathsf T_\rho\mathbf 1_A\rangle
 \ge b^{2/(1-|\rho|)}.
\end{equation}
\end{lemma}

\begin{proof}
First suppose $0\le\rho<1$, so $p=1-\rho$.
We prove the assertion for one coordinate and then pass to products.

For functions on $\{0,1\}$, write $f = a\chi_\emptyset+b\chi_{[1]}$ and $g = c\chi_\emptyset+e\chi_{[1]}$, i.e.,
\begin{equation*}
 f(0)=a+b,\quad f(1)=a-b,
 \qquad
 g(0)=c+e,\quad g(1)=c-e.
\end{equation*}
Nonnegativity gives $a\ge|b|$ and $c\ge|e|$.
Lemma~\ref{appfr:scalar}, applied with $t=|b|/a$ when $a>0$, gives
\begin{equation*}
 \|f\|_p^2\le a^2-\rho b^2.
\end{equation*}
The same inequality is immediate when $a=0$, since then $f=0$.
Likewise $\|g\|_p^2\le c^2-\rho e^2$.
On the other hand,
\begin{equation*}
 \langle f,\mathsf T_\rho g\rangle
 =ac+\rho be\ge ac-\rho|be|.
\end{equation*}
The last expression is nonnegative, and
\begin{equation*}
 (ac-\rho|be|)^2
 -(a^2-\rho b^2)(c^2-\rho e^2)
 =\rho(a|e|-c|b|)^2\ge0.
\end{equation*}
Taking square roots proves the one-coordinate assertion.

Now assume the assertion holds in dimension $n-1$.
Let $(X,Y)$ have the noise coupling with parameter $\rho$: $X$ is uniform, and $Y=X+Z$ with $Z$ distributed as in the definition of $\mathsf T_\rho$.
The coordinate pairs $(X_i,Y_i)$ are independent.
For $i,j\in\{0,1\}$, define the sections
\begin{equation*}
 f_i(x')=f(i,x'),\qquad g_j(y')=g(j,y'),
\end{equation*}
where $x',y'\in\Omega_{n-1}$, and put
\begin{equation*}
 F(i)=\|f_i\|_p,\qquad G(j)=\|g_j\|_p.
\end{equation*}
Conditioning on the first coordinate pair and using the induction hypothesis gives
\begin{alignat*}{2}
 \mathbb E[f(X)g(Y)]
 &\ge\sum_{i,j\in\{0,1\}}
       \mathbb P(X_1=i,Y_1=j)F(i)G(j)\\
 &=\mathbb E[F(X_1)G(Y_1)]\\
 &\ge\|F\|_p\|G\|_p
  =\|f\|_p\|g\|_p.
\end{alignat*}
The second inequality is the one-coordinate assertion.
The last equality follows by expanding the definitions of the $p$th powers.
This proves \eqref{appfr:reverse-bound} for nonnegative $\rho$ in every dimension.

For $\rho=-r<0$, put $\mathbf 1=(1,\ldots,1)$ and $g^\ast(x)=g(x+\mathbf 1)$.
Complementing every coordinate shows directly that
\begin{equation*}
 \mathsf T_{-r}g=\mathsf T_r g^\ast,
 \qquad
 \|g^\ast\|_p=\|g\|_p.
\end{equation*}
Apply the already proved inequality to $f$ and $g^\ast$.
Finally, $\|\mathbf 1_A\|_p=b^{1/p}$, which proves \eqref{appfr:noise-set-bound}.
\end{proof}

\subsection{Elementary estimates for the spherical eigenvalues}

\begin{lemma}
\label{appfr:spectral-estimates}
Let $1\le d\le n-1$, and put
\begin{equation*}
 \delta=\frac dn,\qquad
 \rho_d=1-2\delta,\qquad
 q_d=\delta^2+(1-\delta)^2.
\end{equation*}
For every integer $0\le k\le n$,
\begin{alignat}{2}
 \bigl|\lambda_d(k)-\rho_d^k\bigr|
 &\le\frac{k(k-1)}n,
 \label{appfr:low-degree}\\
 \lambda_d(n-k)&=(-1)^d\lambda_d(k),
 \label{appfr:reflection}\\
 |\lambda_d(k)|&\le(n+1)q_d^{\min(k,n-k)}.
 \label{appfr:middle-degree}
\end{alignat}
In \eqref{appfr:low-degree}, we use the convention $0^0 = 1$ when $d = \frac{n}{2}$.
\end{lemma}

\begin{proof}
\emph{Sampling with and without replacement.}
Consider a list of $n$ signs $\xi_1,\ldots,\xi_n\in\{-1,1\}$, of which
exactly $d$ equal $-1$ and $n-d$ equal $1$.
Let $I_1,\ldots,I_k$ be indices sampled
uniformly without replacement from $[n]$.
If $J$ denotes the number of negative signs
among the sampled values, then
\[
\Prob(J=j)
=
\frac{\binom dj\binom{n-d}{k-j}}{\binom nk}
=
\frac{\binom kj\binom{n-k}{d-j}}{\binom nd}.
\]
The second equality is a standard
binomial-coefficient identity, which can also
be seen by counting pairs of subsets of $[n]$
of cardinalities $d$ and $k$ whose intersection
has cardinality $j$.
Since the product of the sampled signs is
$(-1)^J$, it follows from \eqref{appfr:lambda}
that
\[
\E\left[\prod_{i=1}^k\xi_{I_i}\right]
=
\binom nd^{-1}
\sum_{j=0}^k
(-1)^j\binom kj\binom{n-k}{d-j}
=
\lambda_d(k).
\]

Now let $I_1',\ldots,I_k'$ be independent
indices sampled uniformly from $[n]$, this
time with replacement.
Each sampled sign has expectation
\[
\E[\xi_{I_i'}]
=
\frac{n-d}{n}-\frac dn
=
1-\frac{2d}{n}
=
\rho_d.
\]
By independence,
\[
\E\left[\prod_{i=1}^k\xi_{I_i'}\right]
=
\rho_d^k.
\]

For $k=0$ or $k=1$, the two sampling
procedures give the same expectation, so
\eqref{appfr:low-degree} is immediate.
Suppose therefore that $2\le k\le n$.
Let $D$ be the event that the indices
$I_1',\ldots,I_k'$ are all distinct.
Conditional on $D$, their joint distribution
is exactly that of a uniform ordered sample
without replacement.
Consequently, writing
\[
Y=\prod_{i=1}^k\xi_{I_i'},
\]
we have
\[
\E[Y\mid D]=\lambda_d(k),
\qquad
\E[Y]=\rho_d^k.
\]
By the law of total expectation,
\[
\E[Y]
=
\Prob(D)\E[Y\mid D]
+
\Prob(D^c)\E[Y\mid D^c].
\]
Subtracting $\E[Y\mid D]$ and taking
absolute values gives
\[
\begin{aligned}
|\rho_d^k-\lambda_d(k)|
&=
\Prob(D^c)
\left|
\E[Y\mid D^c]-\E[Y\mid D]
\right|\\
&\le 2\Prob(D^c),
\end{aligned}
\]
where the last inequality follows from
$|Y|=1$.

For each pair of distinct sampling positions
$i<j$, independence and uniformity give
\[
\Prob(I_i'=I_j')=\frac1n.
\]
The event $D^c$ occurs precisely when at
least one such pair of indices coincides.
Thus the union bound yields
\[
\Prob(D^c)
\le
\sum_{1\le i<j\le k}
\Prob(I_i'=I_j')
=
\frac{\binom k2}{n}.
\]
Combining the preceding inequalities, we obtain
\[
|\lambda_d(k)-\rho_d^k|
\le
2\frac{\binom k2}{n}
=
\frac{k(k-1)}n,
\]
which proves \eqref{appfr:low-degree}.

\emph{Reflection of the degree.}
For every $z\in\Omega_n$ with $|z|=d$ and $S \subseteq [n]$,
the definition of the Walsh characters gives
\[
\begin{aligned}
\chi_{S^c}(z)
&=
(-1)^{\sum_{i\in S^c}z_i} =
(-1)^{\sum_{i=1}^n z_i}
(-1)^{\sum_{i\in S}z_i} =
(-1)^d\chi_S(z).
\end{aligned}
\]
By the definition of the spherical
averaging operator and its eigenvalues, we see that for any $S \subseteq [n]$
\[
\lambda_d(|S|)
=
\binom nd^{-1}
\sum_{\substack{z\in\Omega_n\\|z|=d}}
\chi_S(z).
\]
Applying these formulas with $S := [k]$ and $S^c$ gives
\[
\begin{aligned}
\lambda_d(n-k)
&=
\binom nd^{-1}
\sum_{\substack{z\in\Omega_n\\|z|=d}}
\chi_{S^c}(z) =
(-1)^d
\binom nd^{-1}
\sum_{\substack{z\in\Omega_n\\|z|=d}}
\chi_S(z) =
(-1)^d\lambda_d(k),
\end{aligned}
\]
proving \eqref{appfr:reflection}.

\emph{Decay away from the two endpoint degrees.}
By \eqref{appfr:reflection}, it suffices to consider $k\le n/2$.
For a polynomial $P$, write $[z^d]P$ for the coefficient of $z^d$.
Formula~\eqref{appfr:lambda} gives
\begin{equation*}
 \binom nd\lambda_d(k)
 =[z^d](1-z)^k(1+z)^{n-k}
 =[z^d](1-z^2)^k(1+z)^{n-2k}.
\end{equation*}
Taking absolute values term by term yields
\begin{equation*}
 \binom nd|\lambda_d(k)|
 \le[z^d](1+z^2)^k(1+z)^{n-2k}.
\end{equation*}
If a polynomial has nonnegative coefficients, its $d$th coefficient is at most its value at $t$ divided by $t^d$, for every $t>0$.
Therefore
\begin{equation*}
 |\lambda_d(k)|
 \le\frac{t^{-d}(1+t^2)^k(1+t)^{n-2k}}{\binom nd}.
\end{equation*}
Choose $t=\delta/(1-\delta)$.
Since
\begin{equation*}
 \frac{1+t^2}{(1+t)^2}=q_d,
\end{equation*}
we obtain
\begin{equation}
\label{appfr:binomial-factor}
 |\lambda_d(k)|
 \le\frac{q_d^k}{\binom nd\delta^d(1-\delta)^{n-d}}.
\end{equation}
The denominator is the probability of the value $d$ for a binomial random variable with parameters $n$ and $\delta=d/n$.
We claim that $d$ is the value that is most likely to occur out of all possible values.
In other words, we claim that

\begin{equation*}
    \displaystyle\max_{0 \le j \le n}\binom{n}{j}\delta^j(1-\delta)^{n-j} = \binom nd\delta^d(1-\delta)^{n-d}.
\end{equation*}
To see that this is the case, we observe that
\begin{equation*}
 \frac{\binom{n}{j+1}\delta^{j+1}(1-\delta)^{n-j-1}}{\binom{n}{j}\delta^j(1-\delta)^{n-j}} = \frac{n-j}{j+1}\cdot\frac d{n-d},
\end{equation*}
which is greater than one for $j<d$ and less than one for $j\ge d$.
The largest of $n+1$ nonnegative numbers summing to one is at least $1/(n+1)$.
Thus the denominator in \eqref{appfr:binomial-factor} is at least $1/(n+1)$.
This proves \eqref{appfr:middle-degree} for $k\le n/2$, and reflection proves the remaining cases.
\end{proof}

\subsection{Comparing two consecutive distances with independent noise}

For a linear operator $L$ on functions on $\Omega_n$, write
\begin{equation*}
 \|L\|_2=\sup_{f\ne0}\frac{\|Lf\|_2}{\|f\|_2}.
\end{equation*}
In other words, we use $\|L\|_2$ to denote the operator norm of $L$ when viewed as an operator on $L^2(\Omega_n)$.

\begin{lemma}[Operator comparison]
\label{appfr:comparison}
Fix $0<\eta<1/2$.
There exist $C_\eta>0$ and $N_\eta$ such that, whenever $n\ge N_\eta$ and
\begin{equation*}
 \eta n\le d,\qquad d+1\le(1-\eta)n,
\end{equation*}
we have, with $\rho=1-2d/n$,
\begin{equation}
\label{appfr:operator-bound}
 \left\|\frac{\mathsf S_d+\mathsf S_{d+1}}2
               -\mathsf T_\rho\right\|_2
 \le C_\eta\frac{\log^2(n+1)}n.
\end{equation}
\end{lemma}

\begin{proof}
Put
\begin{equation*}
 q=\eta^2+(1-\eta)^2<1,
 \qquad
 \rho'=1-\frac{2(d+1)}n,
 \qquad
 K=\left\lceil\frac{3\log(n+1)}{-\log q}\right\rceil.
\end{equation*}
The hypotheses imply
\begin{equation*}
 q_d,q_{d+1}\le q,
 \qquad
 |\rho|,|\rho'|\le1-2\eta\le q,
 \qquad
 |\rho-\rho'|=\frac2n.
\end{equation*}
Also $q^K\le(n+1)^{-3}$ and $K=O_\eta(\log(n+1))$.
Increase $N_\eta$ so that $2K<n$ whenever $n\ge N_\eta$.

Set
\begin{equation*}
 \mu(k)=\frac{\lambda_d(k)+\lambda_{d+1}(k)}2.
\end{equation*}
By \eqref{appfr:fourier} and \eqref{appfr:eigenvalues}, the operator norm in \eqref{appfr:operator-bound} is
\begin{equation*}
 \max_{0\le k\le n}|\mu(k)-\rho^k|.
\end{equation*}
Indeed, both operators are diagonal in the same orthonormal basis, and the squared norm of their difference applied to $f$ is
\begin{equation*}
 \sum_{S\subseteq[n]}
   |\mu(|S|)-\rho^{|S|}|^2\widehat f(S)^2.
\end{equation*}
We estimate the eigenvalue difference in three ranges.

\emph{Low degrees: $0\le k\le K$.}
For $u,v\in[-1,1]$, factorization of $u^k-v^k$ gives
\begin{equation*}
 |u^k-v^k|\le k|u-v|.
\end{equation*}
Using \eqref{appfr:low-degree} for $d$ and $d+1$, we obtain
\begin{alignat*}{2}
 |\mu(k)-\rho^k|
 &\le\frac12|\lambda_d(k)-\rho^k|
      +\frac12|\lambda_{d+1}(k)-{\rho'}^k|
      +\frac12|{\rho'}^k-\rho^k|\\
 &\le\frac{k(k-1)}n+\frac kn
  =\frac{k^2}n.
\end{alignat*}
The same calculation gives the useful bound
\begin{equation}
\label{appfr:adjacent-low}
 \frac12|\lambda_d(k)-\lambda_{d+1}(k)|\le\frac{k^2}n.
\end{equation}

\emph{Middle degrees: $K<k<n-K$.}
The estimate \eqref{appfr:middle-degree} gives
\begin{equation*}
 |\mu(k)|\le(n+1)q^K,
 \qquad
 |\rho|^k\le q^K.
\end{equation*}
Consequently
\begin{equation*}
 |\mu(k)-\rho^k|\le(n+2)q^K.
\end{equation*}

\emph{High degrees: $n-K\le k\le n$.}
Write $\ell=n-k$, so $0\le\ell\le K$.
Reflection gives
\begin{equation*}
 \mu(k)=\frac{(-1)^d}{2}
                \bigl(\lambda_d(\ell)-\lambda_{d+1}(\ell)\bigr).
\end{equation*}
The opposite signs here are the reason for averaging consecutive distances.
By \eqref{appfr:adjacent-low},
\begin{equation*}
 |\mu(k)|\le\frac{\ell^2}n\le\frac{K^2}n.
\end{equation*}
Moreover, $k\ge n-K>K$ implies $|\rho|^k\le q^K$.
Thus
\begin{equation*}
 |\mu(k)-\rho^k|\le\frac{K^2}n+q^K.
\end{equation*}

Combining the three ranges, we obtain the explicit bound
\begin{equation}
\label{appfr:explicit-error}
 \left\|\frac{\mathsf S_d+\mathsf S_{d+1}}2
               -\mathsf T_\rho\right\|_2
 \le\frac{K^2}n+(n+2)q^K
 \le\frac{K^2}n+\frac{2}{(n+1)^2}.
\end{equation}
Since $K=O_\eta(\log(n+1))$, this proves \eqref{appfr:operator-bound}.
\end{proof}

\begin{remark}\label{rem:DensityFRErrorBound}
    Combining the exponential Krawtchouk bound of Harrow--Kolla--Schulman \cite[Lemma~2.2, p.~60]{HKS2014} with the classical Krawtchouk recurrence recorded in \cite[equation~(2.26), p.~741]{Polyanskiy2019}, one can show that

    \begin{equation}
\label{appfr:operator-boundInRemark}
 \frac{1}{n} \le \left\|\frac{\mathsf S_d+\mathsf S_{d+1}}2
               -\mathsf T_\rho\right\|_2
 \le \frac{C_\eta}n.
\end{equation}
    This means that the error bound in Equation \eqref{eq:BHM} can be taken to be $\frac{C_\gamma}{n}$ instead of $C_\gamma\frac{\log^2(n+1)}{n}$.
    We do not prove this sharper bound here because it does not help with our applications to recurrence, and the goal of this appendix is to give a self-contained proof of Theorem \ref{thm:BHM}.
    We also mention that the error bound obtained from the original proof of Benabbas \cite[Theorem 7.1]{BenabbasThesis} is $C_\gamma n^{-\frac{1}{5}}$.
\end{remark}

\subsection{Parity reduction and completion of the proof}

\begin{proof}[Proof of Theorem~\ref{thm:BHM}]
Fix $0<\gamma<1/2$ and an $n$ for which $d=(1-\gamma)n$ is an even integer.
Put
\begin{equation*}
 \rho=1-\frac{2d}n=2\gamma-1,
 \qquad
 \mathsf A_d=\frac{\mathsf S_d+\mathsf S_{d+1}}2.
\end{equation*}
Apply Lemma~\ref{appfr:comparison} with $\eta=\gamma/2$.
Its hypotheses hold for all sufficiently large $n$, because $d/n=1-\gamma$ and $1/n\le\gamma/2$ eventually.
It follows that there is some $\varepsilon_{n,\gamma}$ for which
\begin{equation}
\label{appfr:final-error}
 \|\mathsf A_d-\mathsf T_\rho\|_2
 \le \varepsilon_{n,\gamma} \le C_\gamma\frac{\log^2(n+1)}n.
\end{equation}

Now partition $B$ according to the parity of its Hamming weight:
\begin{equation*}
 B=B_0\sqcup B_1,
 \qquad
 B_i=\{x\in B:|x|\equiv i\pmod2\}.
\end{equation*}
Choose $i$ so that $|B_i|\ge|B|/2$, and put
\begin{equation*}
 f=\mathbf 1_{B_i},
 \qquad
 b=\mathbb E f=\frac{|B_i|}{2^n}.
\end{equation*}
We see that $\frac{a}{2} \le b\le \frac{1}{2}$, and that any two points in $B_i$ have even Hamming distance.
Since $d+1$ is odd, we see that
\begin{equation*}
 \langle f,\mathsf S_{d+1}f\rangle=0,
 \qquad
 \langle f,\mathsf S_df\rangle
       =2\langle f,\mathsf A_df\rangle.
\end{equation*}
By Cauchy--Schwarz and \eqref{appfr:final-error},
\begin{equation*}
 \left|\langle f,(\mathsf A_d-\mathsf T_\rho)f\rangle\right|
 \le \varepsilon_{n,\gamma}\|f\|_2^2
 = \varepsilon_{n,\gamma}b.
\end{equation*}
Also $1-|\rho|=2\gamma$, so Lemma~\ref{appfr:reverse} gives
\begin{equation*}
 \langle f,\mathsf T_\rho f\rangle\ge b^{1/\gamma}.
\end{equation*}
Since $B_i\subseteq B$ and the spherical averaging operator has nonnegative coefficients, we conclude that
\begin{alignat*}{2}
 \langle\mathbf 1_B,\mathsf S_d\mathbf 1_B\rangle
 &\ge\langle f,\mathsf S_df\rangle =2\langle f,\mathsf A_df\rangle\\
 &\ge2\langle f,\mathsf T_\rho f\rangle
      -2\varepsilon_{n,\gamma}b\\
 &\ge2b^{1/\gamma}-2\varepsilon_{n,\gamma}b\\
 &\ge2(a/2)^{1/\gamma}-C_\gamma\frac{\log^2(n+1)}n.
\end{alignat*}
Together with \eqref{appfr:edge-probability}, this is \eqref{eq:BHM}.
\end{proof}

\begin{remark}
The assumption that $d$ is even is essential: a full parity class has density $1/2$ but contains no pair at an odd Hamming distance.
\end{remark}

\section{Sharp \texorpdfstring{$s$}{s}-nice recurrence without the vdC property in \texorpdfstring{$\mathbb{F}_2^\omega$}{F2-omega}}
\label{app:boolean-recurrence}

Let $G=\mathbb F_2^\omega=\bigoplus_{j\ge1}\mathbb F_2$.
For a subset of this discrete group, the notation $r\to\infty$ means that $r$ eventually leaves every finite subset of $G$.
We use the natural analogue of $s$-nice measurable recurrence for measure-preserving $G$-actions.
We use the dynamical characterization of vdC sets given in \cite[Theorem 4.2]{FarhangiTuckerDrob2026} that $R \subseteq G$ is a vdC set if for every p.p.s. $(X,\mathscr{B},\mu,(T^g)_{g \in G})$ and every $f \in L^\infty(X,\mu)$ satisfying $\langle f\circ T^r,f\rangle = 0$ for all $r \in R$, we have $\int_Xfd\mu = 0$.\footnote{In \cite{RodriguezMartin2025}, this dynamical characterization is taken as the definition of a vdC set in an arbitrary countably infinite group $G$.}

\begin{theorem}\label{thm:boolean-separation}
For every real $s>2$ there is a set $R_s\subseteq G\setminus\{0\}$ that is $t$-nicely recurrent exactly when $t\ge s$, but is not a vdC set.
More precisely, there are a probability-preserving $G$-action $(T^g)_{g\in G}$ and a real function $f\in L^2(X,\mu)$ such that
\begin{equation}\label{eq:boolean-witness}
 \|f\|_2=1,\qquad
 \int f\,d\mu=\sqrt{\frac{s-2}{2(s-1)}}>0,\qquad
 \langle f\circ T^r,f\rangle=0\quad(r\in R_s).
\end{equation}
\end{theorem}

\begin{proof}
Choose rational numbers $\gamma_j\uparrow1/s$ with $0<\gamma_j<1/s$, and put $p_j=1/\gamma_j$ and $c_j=2^{1-p_j}$.
By Theorem~\ref{thm:BHM}, choose dimensions $n_j\to\infty$ such that $d_j=(1-\gamma_j)n_j$ is even and, for every $B\subseteq\mathbb F_2^{n_j}$ of density $a$,
\begin{equation}\label{eq:boolean-density}
 \frac{\#\{(u,v)\in B^2:|u-v|=d_j\}}
 {2^{n_j}\binom{n_j}{d_j}}
 \ge c_j a^{p_j}-2^{-j}.
\end{equation}
Identify $G$ with $\bigoplus_{j\ge1}G_j$, where $G_j=\mathbb F_2^{n_j}$, and define
\[
 E_j=\{r\in G_j:|r|=d_j\},\qquad R_s=\bigcup_{j\ge1}E_j.
\]
The blocks are disjoint and escape every finite subset of $G$.

\emph{$s$-Nice recurrence.}
Let $(X,\mathscr{B},\mu,(T^g)_{g \in G})$ be a p.p.s., let $A \in \mathscr{B}^+$ be arbitrary, and put $a=\mu(A)>0$ and $M_j(A)=\max_{r\in E_j}\mu(A\cap T^{-r}A)$.
For $x \in X$, set $B_j(x)=\{v\in G_j:T^v x\in A\}$ and $b_j(x)=2^{-n_j}|B_j(x)|$.
Since $T$ preserves $\mu$, we see that $\int b_j\,d\mu=a$.
Averaging \eqref{eq:boolean-density} over $x$ and applying Jensen's inequality yields
\[
 M_j(A)\ge c_j a^{p_j}-2^{-j}.
\]
Applying the same argument to the $k$-fold product action and the event $A^k$ gives, for every integer $k\ge1$,
\begin{equation}\label{eq:boolean-tensor}
 M_j(A)^k\ge c_j a^{kp_j}-2^{-j}.
\end{equation}
Since $p_j\to s$, $c_j\to2^{1-s}$, and $2^{-j}\to0$, it follows that
\[
 \liminf_{j\to\infty}M_j(A)\ge 2^{(1-s)/k}a^s.
\]
Letting $k\to\infty$ shows that $R_s$ is $s$-nicely recurrent.

\emph{non-vdC witness.}
Set $\alpha=1-2/s\in(0,1)$ and $a_j=2d_j/n_j-1=1-2\gamma_j>\alpha$.
For $r\in E_j$, the degree-one Walsh characters $\chi_{j,i}(v)=(-1)^{v_i}$ satisfy the pointwise identity
\begin{equation}\label{eq:boolean-degree-one}
 \frac1{n_j}\sum_{i=1}^{n_j}\chi_{j,i}(r)
 =\frac{n_j-2d_j}{n_j}=-a_j.
\end{equation}
We remark that $-a_j$ is also the eigenvalue of the spherical averaging operator $\mathsf S_{d_j}$ corresponding to the Walsh characters $\chi_{j,i}$.
Let $u_j$ be the uniform distribution on these $n_j$ characters, and let $m_j$ be the normalized Haar measure on the dual group $\widehat G_j$.
The probability measure
\[
 \rho_j=\frac\alpha{a_j}u_j+\left(1-\frac\alpha{a_j}\right)m_j
\]
therefore satisfies $\int_{\widehat{G}_j}\chi(r)d\rho_j(\chi) =-\alpha$ for every $r\in E_j$, because the average of all characters at a nonzero group element is zero.
For $\omega \in \Omega := \prod_{j \ge 1}\widehat G_j$ and  $g=(g_j)\in G$, let $\chi_\omega(g)=\prod_{j \ge 1}\omega_j(g_j)$.
Give $X=\Omega\times\{-1,1\}$ the product of $\rho := \otimes_{j \ge 1}\rho_j$ with the uniform measure on $\{-1,1\}$, and put
\[
 T^g(\omega,\varepsilon)=(\omega,\chi_\omega(g)\varepsilon),
 \qquad h(\omega,\varepsilon)=\varepsilon.
\]
These maps form a measure-preserving $G$-action, and $\int h\,d\mu=0$, $\|h\|_2=1$, and
\[
 \langle h\circ T^r,h\rangle
 =\int_\Omega\chi_\omega(r)d\rho(\omega)=-\alpha\qquad(r\in R_s).
\]
Thus, the function
\begin{equation}\label{eq:boolean-f}
 f=\frac{\sqrt\alpha+h}{\sqrt{1+\alpha}}
\end{equation}
satisfies $\|f\|_2=1$, $\int f\,d\mu=\sqrt{\alpha/(1+\alpha)}$, and
\[
 \langle f\circ T^r,f\rangle
 =\frac{\alpha+\langle h\circ T^r,h\rangle}{1+\alpha}=0
 \qquad(r\in R_s).
\]
Its nonzero mean witnesses failure of the vdC property.

\emph{Sharpness of the exponent.}
The Gaussian measure space construction that we discussed in the proof of Proposition~\ref{prop:gaussian-obstruction} applies when $\mathbb{Z}$ is replaced by any other countable group, such as $G$.
We consider the real positive definite function $c:G\rightarrow\mathbb{R}$ given by $c(g)=\langle h\circ T^g,h\rangle$, and we note that $c(0)=1$ and $c(r)=-\alpha$ for $r\in R_s$.
Let $(Y_g)_{g\in \mathbb{F}_2^\omega}$ be centered Gaussian variables with covariance $\mathbb E(Y_gY_h)=c(g-h)$, and let $A_u=\{Y_0>u\}$.
For $r\in R_s$, the variance of $Y_0+Y_r$ is $2(1-\alpha)=4/s$, so the Gaussian exponential bounds from before give
\[
 \mu(A_u\cap T^{-r}A_u)
 \le \mathbb P(Y_0+Y_r>2u)\le\exp(-su^2/2).
\]
Since $\mu(A_u)\ge C u^{-1}e^{-u^2/2}$ for $u\ge1$ and some absolute $C>0$, we obtain
\[
 \sup_{r\in R_s}\frac{\mu(A_u\cap T^{-r}A_u)}{\mu(A_u)^t}
 \le C^{-t}u^t\exp\!\left(-\frac{(s-t)u^2}{2}\right)
 \longrightarrow0\qquad(u\to\infty)
\]
for every $0<t<s$.
In particular, $R_s$ is not $t$-nicely recurrent for any $2\le t<s$.
\end{proof}

\begin{remark}[The CRT transfer from Hamming cubes to integers]
\label{rem:TransferFromTorsion}
The strategy of constructing examples in finite high-dimensional groups and then transferring them to $\mathbb{Z}$ has several precedents in the literature on recurrence.
A particularly close precedent for our CRT construction is that of Ruzsa \cite[Sections~6--7]{Ruzsa1985Bohr}, who uses distinct odd prime moduli to encode finite combinatorial configurations with coordinatewise differences taking values in $\{0,\pm1\}$.
The broader strategy of transferring finite high-dimensional configurations to $\mathbb{Z}$ also appears in the work of K\v{r}\'{\i}\v{z} \cite{Kriz1987}; see \cite[Sections~1.2, 4, and~5]{Griesmer2025Topological} for an exposition.
Related transfer arguments involving the approximation of characters and the preservation of densities can be found in \cite[Section~7]{Griesmer2019Rigidity} and \cite[Sections~6--7]{Griesmer2021Bohr}, respectively.
The same general idea, implemented through various encoding and approximation techniques, is also present in the constructions of Bourgain, Mountakis, McCutcheon, and Forrest \cite{Bourgain,Mountakis,McCutcheon1995,Forrest}.

In our setting, the CRT blocks construction allows us to transfer the natural construction in $\mathbb{F}_2^\omega$ to $\mathbb{Z}$ by approximately realizing the relevant Fourier-theoretic properties of each finite Boolean group $G_j$.
Indeed, suppressing the block index $j$, for distinct odd primes $q_1,\ldots,q_n$ and $Q=\prod_{i=1}^n q_i$, we define the injective map $\Phi:\mathbb{F}_2^n\rightarrow\{0,\ldots,Q-1\}$ by requiring that $\Phi(v)\equiv v_i\pmod{q_i}$ for every $1\leq i\leq n$.
For $u,v\in\mathbb{F}_2^n$ and $r=\Phi(v)-\Phi(u)$, the residue $r\bmod q_i$ belongs to $\{0,\pm1\}$ and is nonzero precisely when $u_i\neq v_i$.
Thus, although $\Phi$ is not a group homomorphism, it preserves the support pattern of coordinatewise differences and allows us to recover the corresponding Hamming distances from the CRT coordinates.
Moreover, if $t_i=(q_i-1)/(2q_i)$, then
\[
\cos(2\pi rt_i)=
\begin{cases}
    1,&u_i=v_i,\\
    -\cos(\pi/q_i),&u_i\neq v_i,
\end{cases}
\]
and consequently
\[
\left|\cos(2\pi rt_i)-(-1)^{u_i+v_i}\right|
\leq\frac{\pi^2}{2q_i^2}.
\]
Recalling that $\chi_i(v-u)=(-1)^{u_i+v_i}$, we see that the Fourier coefficients of the point masses $\delta_{\chi_i}\in\mathcal{M}(\widehat{G}_j)$ are approximated, along the corresponding encoded differences, by those of the symmetric measures $\frac{1}{2}(\delta_{t_i}+\delta_{-t_i})\in\mathcal{M}(\mathbb{T})$, with an error that can be made arbitrarily small by choosing the primes $q_i$ sufficiently large.
It is in this Fourier-theoretic sense that the CRT blocks construction provides an approximate realization of the finite Boolean models in $\mathbb{Z}$.

While the transfer of the Fourier-theoretic properties is approximate, the transfer of the combinatorial properties responsible for $s$-nice recurrence is exact.
Indeed, writing $E_j=\{g\in G_j:|g|=d_j\}$, we see that $\Phi_j$ is an injective graph homomorphism from $\operatorname{Cay}(G_j,E_j)$ to $\operatorname{Cay}(\mathbb{Z},D_j\cup(-D_j))$.
This preservation of the relevant combinatorial configurations allows the same orbit-averaging and Cartesian-power arguments used to establish $s$-nice recurrence in $\mathbb{F}_2^\omega$ to establish $s$-nice recurrence in $\mathbb{Z}$.
\end{remark}

\begin{remark}\label{rem:OptimalAtom}
    We mention that Proposition~\ref{prop:gaussian-obstruction} holds when $\mathbb{Z}$ is replaced by any countable group $G$.
    Consequently, we can use Proposition~\ref{prop:gaussian-obstruction} for $\mathbb{F}_2^\omega$ to prove the analogue of Theorem~\ref{thm:OptimalAtom} in $\mathbb{F}_2^\omega$, which can be used to show that the value of $\int fd\mu$ appearing in Equation \eqref{eq:boolean-witness} is the largest possible.
    The only adjustment we need to make to the proof of Theorem~\ref{thm:OptimalAtom} is how we handle the case of $\beta = \frac{1}{2}$.
    We see that if $\beta = \sigma(\{0\}) = \frac{1}{2}$ and $\widehat{\sigma}(r_0) = 0$, then for $\nu = \sigma-\frac{1}{2}\delta_0$ we have $\text{supp}(\nu) \subseteq \{\chi \in \widehat{\mathbb{F}_2^{\omega}}\ |\ \chi(r_0) = -1\}$.
    Let $\chi_0 \in \text{supp}(\nu)$ be arbitrary, and consider the probability-preserving action $T$ of $\mathbb{F}_2^\omega$ on $\{-1,1\}$ given by $T^g(i) = \chi_0(g)i$.
    Since $R$ is a set of recurrence, there exists $r_1 \in R$ for which $\chi_0(r_1) = 1$.
    Since the evaluation map $\chi\mapsto\chi(r)$ is continuous, we see that there is an open neighborhood $U$ of $\chi_0$ such that $\chi(r_1) = 1$ for every $\chi \in U$, which shows us that $\widehat{\sigma}(r_1) > 0$.
    This yields the desired contradiction and shows that $\beta < \frac{1}{2}$.
    \end{remark}

\end{document}